\documentclass[11pt]{amsart}
\usepackage{amssymb}
\usepackage{epsfig}
\usepackage{amsfonts}
\usepackage{indentfirst}
\usepackage{epsf}

\numberwithin{equation}{section}
\newtheorem{theorem}{Theorem}[section]
\newtheorem{lemma}[theorem]{Lemma}

\newtheorem{corollary}[theorem]{Corollary}

\theoremstyle{definition}
\newtheorem{definition}[theorem]{Definition}
\newtheorem{remark}[theorem]{Remark}
\newtheorem{example}[theorem]{Example}

\DeclareMathOperator{\pd}{pd}
\DeclareMathOperator{\ara}{ara}
\DeclareMathOperator{\height}{height}
\DeclareMathOperator{\level}{level}

\DeclareMathOperator{\dist}{dist}

\begin{document}

\title{Algebraic properties of the path ideal of a tree}

\author{Jing (Jane) He}
\address{School of Mathematics \& Statistics\\
Carleton University\\
Ottawa, ON K1S 5B6, Canada}
\email{jhe2@connect.carleton.ca}

\author{Adam Van Tuyl}
\address{Department of Mathematical Sciences \\
Lakehead University \\
Thunder Bay, ON P7B 5E1, Canada}
\email{avantuyl@lakeheadu.ca}
\urladdr{http://flash.lakeheadu.ca/$\sim$avantuyl/}

\thanks{Revised Version: April 20, 2009}

\begin{abstract}
The path ideal (of length $t \geq 2$) of a directed graph $\Gamma$
is the monomial ideal, denoted $I_t(\Gamma)$, whose generators correspond
to the directed paths of length $t$ in $\Gamma$.  We study
some of the algebraic properties of $I_t(\Gamma)$ when $\Gamma$ is a tree.
We first show that $I_t(\Gamma)$ is the facet ideal of a simplicial tree.
As a consequence, the quotient ring $R/I_t(\Gamma)$ is always sequentially
Cohen-Macaulay, and the Betti numbers of $R/I_t(\Gamma)$ do not depend
upon the characteristic of the field.
We study the case of the line graph in greater detail at the
end of the paper.
We give an exact formula for the projective dimension of these
ideals, and in some cases,
we compute their arithmetical rank.
\end{abstract}

\maketitle


\section{Introduction}

Let $\Gamma = (V,E)$ be a finite simple graph with vertex set
$V = \{x_1,\ldots,
x_n\}$ and edge set $E$.  By identifying the vertices with the variables
in the polynomial ring $R = k[x_1,\ldots,x_n]$ over a field $k$,
there exists a growing number of ways to associate to a graph $\Gamma$ a
(monomial) ideal in $R$.  The best known correspondence
is the edge ideal construction of Villarreal \cite{V1}
where $G$ is associated to the monomial ideal whose
generators correspond to the edges of $\Gamma$.
Other constructions and higher dimensional analogues can
be found in \cite{CFF1,F1,HVT1}, among others.  The underlying theme
in all correspondences is to relate the algebraic properties
to the combinatorial properties, and vice versa.
We contribute to this program by studying the algebraic
properties of the path ideal.

The path ideal of a graph was first introduced by Conca
and De Negri \cite{CD}.
Fix an integer $t\geq 2$, and suppose that
$\Gamma$ is a directed graph, i.e., each edge
has been assigned a direction.  A sequence of $t$ vertices
$x_{i_1},\ldots,x_{i_t}$ is said to be a {\bf path} of length
$t$ if there are $t-1$ distinct edges $e_1,\ldots,e_{t-1}$ such that
$e_j = (x_{i_j},x_{i_{j+1}})$ is a directed
edge from $x_{i_j}$ to $x_{i_{j+1}}$.  The {\bf path ideal}
of $\Gamma$ of length $t$ is the monomial ideal
\[I_t(\Gamma) = \left(\{x_{i_1}\cdots x_{i_t} ~|~ x_{i_1},\ldots, x_{i_t}
~\mbox{is a path of length $t$ in $\Gamma$}\}\right).\]
Note that when $t=2$, then $I_2(\Gamma)$ is simply the
edge ideal of $\Gamma$, and thus $I_t(\Gamma)$ is sometimes
called the {\bf generalized edge ideal} of $\Gamma$.

Path ideals appeared in \cite{CD} as an example
of a family of monomial ideals that are generated by $M$-sequences.
Among other things, it is shown that when $\Gamma$ is a directed tree,
the Rees algebra $\mathcal{R}(I_t(\Gamma))$
is normal and Cohen-Macaulay.  The path ideals of complete bipartite
graphs are shown to be normal in \cite{RV1}, while
the path ideals of cycles are shown to have linear type
in \cite{BdS1}.   To the best of our knowledge, little else
is known about these ideals.  It is therefore of interest
to determine further algebraic properties of the ideals $I_t(\Gamma)$.
In this paper we shall focus on the case that $\Gamma$ is a directed tree,
where the directions of the edges will depend upon which vertex
is designated the root.
By restricting to trees, we can exploit the fact
that there is a unique path between any two vertices.  As we
shall show, many of the results known to hold when $t =2$, that is,
the edge ideal of a tree, continue to hold when $t > 2$.

When $\Gamma$ is a tree, the ideal $I_t(\Gamma)$ is a square-free
monomial ideal.  We can then view $I_t(\Gamma)$ as the facet
ideal of a pure simplicial complex $\Delta_t(\Gamma)$ where
the facets of $\Delta_t(\Gamma)$ correspond to the paths of length
$t$ in $\Gamma$.   Our first main result (see Theorem \ref{simplicialtree})
is that $\Delta_t(\Gamma)$ is a simplicial tree, as defined by
Faridi \cite{F1}.   In other words, $I_t(\Gamma)$ is the facet
ideal of a simplicial tree.

Once we have established that $I_t(\Gamma)$ is the facet ideal
of a simplicial tree, we can employ some well known results
about these ideals.  For example, the main
result of Faridi \cite{F2} implies:

\begin{theorem}[Corollary \ref{scm}]
If $I_t(\Gamma)$ with $t \geq 2$ is the path ideal of a tree $\Gamma$, then
$R/I_t(\Gamma)$ is sequentially Cohen-Macaulay.
\end{theorem}

\noindent
Similarly, results about simplicial forests due to the second
author and H.T. H\`a \cite{HVT2} imply that the graded Betti
numbers of $I_t(\Gamma)$ are independent of the characteristic
of the field $k$ (see Theorem \ref{characteristic}).  Furthermore, the notion of a Betti-splitting,
as defined in \cite{FHVT1}, is used to derive a recursive formula for the
projective dimension for a special subclass of these ideals
(Theorem \ref{generalprojdim}).  Our subclass contains
the edge ideals of trees, and thus, our formula generalizes
a result of Jacques and Katzman \cite{JK1}.

The last part of the paper focuses on the case that $\Gamma = L_n$
is the line graph of length $n$.  In this situation, we derive an exact formula
for the projective dimension of $I_t(L_n)$:

\begin{theorem}[Theorem \ref{maintheorem1}]
Fix integers $n \geq t \geq 2$, and let $L_n$ be the line graph. Then
\[\operatorname{pd}(R/I_t(L_n))=
\left\{
\begin{array}{ll}
\displaystyle{\frac{2(n-d)}{t+1}} & \hbox{if $n\equiv d \pmod{(t+1)}$ with $0 \leq d \leq t-1$} \\
\displaystyle{\frac{2n - (t-1)}{t+1}} & \hbox{if $n\equiv t
\pmod{(t+1)}$.}
\end{array}
\right.
\]
\end{theorem}

\noindent
By specializing this result to the case that $t = 2$, we recover
a result of Jacques \cite{J1}.  We end this paper by computing
the arithmetical rank of $I_t(L_n)$ for {\it some} $t$ and $n$.
Our work extends the work of Barile \cite{B1} who computed
the arithmetical rank of $I_2(L_n)$ for all $n$.  We apply these
results by giving a lower bound (Theorem \ref{lowerbound}) on $\pd(I_t(\Gamma))$ for any
tree $\Gamma$.


\section{Path ideals as facet ideals of simplicial trees}

In this section we show that when $\Gamma$ is a tree,
the ideal $I_t(\Gamma)$ can be viewed as the facet ideal
of a simplicial tree.  Recall that a graph $\Gamma$ is a {\bf tree}
if for every pair of distinct vertices $x$ and $y$,
there exists a unique path in $\Gamma$ between $x$ and $y$.

Throughout this paper, we make the convention
that all of our trees $\Gamma$ are directed graphs.  We recall
the relevant definitions:

\begin{definition} A {\bf directed edge} of a graph
is an assignment of a direction to an edge of a graph.
If $\{w,u\}$ is an edge, we write $(w,u)$ to denote
the directed edge where the direction is
from $w$ to $u$.
A graph is a {\bf directed graph} if each edge has been
assigned a direction.
A {\bf path} of length $t$ in a directed graph
is a sequence of vertices $x_{i_1},\ldots,x_{i_t}$
such that $e_j = (x_{i_j},x_{i_j+1})$ is a directed edge
for $j = 1,\ldots,t-1$.
\end{definition}

A tree $\Gamma$ can be viewed as a directed graph by
picking any vertex of $\Gamma$ to be the {\bf root}
of the tree, and assigning to each edge the direction
``away'' from the root.  Because $\Gamma$ is a tree,
the assignment of a direction will always be possible.
A {\bf leaf} is any vertex in $\Gamma$ adjacent
to only one other vertex.   The {\bf level} of a vertex $x$,
denoted $\level(x)$,
is length of the unique path starting at the root and
ending at $x$ minus one.  The {\bf height} of a tree,
denoted $\height(\Gamma)$, is then given by
$\height(\Gamma) := \max_{x \in V} \{\level(x)\}$.

\begin{example}\label{treeexample}
Consider the following tree $\Gamma$:
\[
\begin{picture}(180,180)(0,0)
\put(100,160){\circle*{5} $x_1$}
\put(100,160){\line(1,-1){40}}
\put(100,160){\vector(1,-1){20}}
\put(100,160){\line(-1,-2){60}}
\put(100,160){\vector(-1,-2){30}}
\put(100,160){\vector(-1,-2){10}}
\put(100,160){\vector(-1,-2){50}}
\put(80,120){\circle*{5} $x_2$}
\put(60,80){\circle*{5} $x_4$}
\put(60,80){\line(1,-2){20}}
\put(60,80){\vector(1,-2){10}}
\put(40,40){\circle*{5} $x_8$}
\put(80,40){\line(-1,-2){20}}
\put(80,40){\vector(-1,-2){10}}
\put(80,40){\circle*{5} $x_9$}
\put(60,0){\circle*{5} $x_{12}$}
\put(120,80){\circle*{5} $x_5$}
\put(140,120){\circle*{5} $x_3$}
\put(140,120){\line(-1,-2){20}}
\put(140,120){\vector(-1,-2){10}}
\put(140,80){\circle*{5}$x_6$}
\put(140,120){\line(0,-1){80}}
\put(140,120){\vector(0,-1){20}}
\put(140,120){\vector(0,-1){60}}
\put(140,120){\line(1,-2){20}}
\put(140,120){\vector(1,-2){10}}
\put(140,80){\line(1,-1){40}}
\put(140,80){\vector(1,-1){20}}
\put(160,80){\circle*{5} $x_7$}
\put(180,40){\circle*{5} $x_{11}$}
\put(140,40){\circle*{5} $x_{10}$}
\end{picture}
\]
The tree $\Gamma$ is a rooted tree whose root is $x_1$.
Every edge has been assigned a direction ``away'' from the root,
as denoted by the arrows.
If $t = 3$, then the path ideal of $\Gamma$ is given by
\[I_3(\Gamma) = (x_1x_2x_4,x_2x_4x_8,x_2x_4x_9,x_4x_9x_{12},x_1x_3x_5,x_1x_3x_6,x_1x_3x_7,
x_3x_6x_{10},x_3x_6x_{11}).\]
The leaves of $\Gamma$ are the vertices: $x_5,x_7,x_8, x_{10},x_{11},x_{12}$.
Observe that $\level(x_{12}) =4$
since the length of the path from $x_1$ to $x_{12}$ has length 5.
By observation, $\height(\Gamma) = \level(x_{12}) = 4$.

Note that
although $x_2,x_1,x_3$ is a path in the {\it undirected} graph $\Gamma$,
we do not consider $x_2x_1x_3$ as a generator of $I_3(\Gamma)$.
If we picked another vertex to be a root, say $x_4$, then the path
ideal is different; in particular:
\[I_3(\Gamma) = (x_4x_9x_{12},x_4x_2x_1,x_2x_1x_3,x_1x_3x_5,x_1x_3x_6,x_1x_3x_7,x_3x_6x_{10},x_3x_6x_{11}).\]
When $x_4$ is the root, the height of $\Gamma$ also changes.  In this context,
$\height(\Gamma) = 5$.
\end{example}

A {\bf simplicial complex} $\Delta$ on the vertex set $V =
\{x_1,\ldots,x_n\}$ is a collection of subsets of $V$ such that:
$(i)$ if $F \in \Delta$ and $G \subseteq F$, then $G \in \Delta$,
and $(ii)$ $\{x_i\} \in \Delta$ for $i = 1,\ldots,n$.  An element $F
\in \Delta$ is a {\bf facet} if $F$ is maximal with respect to
inclusion.  If $\{F_1,\ldots,F_l\}$ is a complete list of the facets
of $\Delta$, then we usually write $\Delta = \langle
F_1,\ldots,F_l\rangle$.  Faridi \cite{F1} introduced the notion of a
{\bf facet ideal} of a simplicial complex; precisely:
\[I(\Delta) := \left.\left(\left\{ \prod_{x \in F} x ~\right|~ \mbox{$F$ a facet of $\Delta$}\right\}\right).\]
Using the facet ideal correspondence, there is a bijection
between the square-free monomial ideals and simplicial
complexes.

Because $I_t(\Gamma)$ is a square-free monomial ideal, this
ideal corresponds to a simplicial complex, say $\Delta_t(\Gamma)$.
The facets of $\Delta_t(\Gamma)$
are precisely the paths of length $t$ in $\Gamma$.  That is,
\[\Delta_t(\Gamma) :=
\left.\left\langle \left\{ \{x_{i_1},\ldots,x_{i_t}\} ~\right|~
\mbox{$x_{i_1},\ldots,x_{i_t}$ is a path of length $t$ in $\Gamma$ }
\right\}\right\rangle.\]
The definitions immediately imply that $I_t(\Gamma) = I(\Delta_t(\Gamma))$.
Our aim is to show that $\Delta_t(\Gamma)$ is a simplicial tree,
as first defined by Faridi \cite{F1}.

\begin{definition}\label{simplicial tree}
A {\bf leaf} of a simplicial complex $\Delta$ is a facet $F$
of $\Delta$ such that either $F$ is the only facet of $\Delta$, or there
exists a facet
$G$ in $\Delta$, $G \not= F$, such that $F \cap F' \subseteq F \cap G$
for every
facet $F' \in \Delta, F' \not= F$.  A simplicial complex $\Delta$
is a {\bf simplicial tree} if $\Delta$ is connected, and every
nonempty subcomplex $\Delta'$
contains a leaf.  By a subcomplex, we mean any simplicial
complex of the form $\Delta' = \langle F_{i_1},\ldots,F_{i_q}\rangle$
where $F_{i_1},\ldots,F_{i_q}$ is a subset of the facets of $\Delta$.
We adapt the convention that the empty simplicial complex
is also a tree.
\end{definition}

\begin{remark}
Note that the term {\bf leaf} will be used in two different
ways.  When referring to a tree $\Gamma$, a leaf is a vertex.  However,
when referring to a simplicial complex, a leaf is a facet.
\end{remark}

We begin with two lemmas about paths in $\Gamma$.  The first lemma
proves that the intersection of any two paths
must start at the first vertex of one of these two paths.

\begin{lemma}\label{intersectionoftwopaths} The intersection of any two
distinct paths $F$ and $G$
in a directed rooted tree $\Gamma$ is a
connected path of length $|F\cap G|$ starting at the first vertex of either
$F$ or $G$.  Furthermore,  the first vertex of this new path
is the first vertex of $F$ and $G$ whose level is the largest.
\end {lemma}

\begin{proof}
Let $F = x_{i_1},\ldots,x_{i_t}$ denote the first path,
and let $G = x_{j_1},\ldots,x_{j_r}$ denote the second path.
We abuse notation and write $F \cap G$ to denote the vertices
that appear in both $F$ and $G$.  Note that if $F \cap G = \emptyset$,
then there is nothing to prove.

Assume $F \cap
G$ starts at $x_m$, which is neither the first vertex of $F$ nor $G$.
So, there are vertices $u\in F$ and $w\in G$ such that $(u,x_m)$ is a
directed edge in the path $F$ and $(w,x_m)$ is a directed edge in $G$.
Because $\Gamma$ is a tree, there is a distinct
path from $x$ to $u$ and also one from $x$ to $w$, where
$x$ is the root.  But then there
are two distinct paths from $x$ to $x_m$, that is, one that goes
through $u$, and the other that goes through $w$.  This
contradicts the definition of a tree. So $F\cap G $ starts at the
first vertex of $F$ or $G$.  It is straightforward
to now see that the first vertex of $F \cap G$ is
$x_k$ where $k = \max\{i_1,j_1\}$.

Suppose $F \cap G = x_{i_{m_1}},\ldots, x_{i_{m_s}}$.
We claim that $F \cap G$ is a
connected path in $\Gamma$. Suppose $F \cap G$ is
the union of two disjoint paths, that is,  $x_{i_{m_1}},\ldots,
x_{i_{m_j}} \cup x_{i_{m_{j+1}}},\ldots, x_{i_{m_s}}$, where
$x_{i_{m_1}},\ldots, x_{i_{m_j}}$ is a path, and
$x_{i_{m_{j+1}}},\ldots, x_{i_{m_s}}$
is path, but there is no edge $(x_{i_{m_j}},x_{i_{m_{j+1}}})$.
Because $F$ is a connected path, there exists vertices
$x_{g_1},\ldots,x_{g_a}$ in $F$ such that
there is path $ x_{i_{m_j}},x_{g_1},\ldots,x_{g_a},x_{i_{m_{j+1}}}$
in $F$.
Similarly, there are vertices $x_{h_1},\ldots,x_{h_b}$ in $G$
that give a path $x_{i_{m_j}},x_{h_1},\ldots,
x_{h_b},x_{i_{m_{j+1}}}$ in $G$.  Note that $x_{g_1},\ldots,x_{g_a},
x_{h_1},\ldots,x_{h_b}$ cannot all be in $F \cap G$, because
if they were, then there would be a path between $x_{i_{m_j}}$
and $x_{i_{m_{j+1}}}$ in $F \cap G$.
But then there are two different paths
between $x_{i_{m_j}}$ and $x_{i_{m_{j+1}}}$ in $\Gamma$, contradicting
the fact that $\Gamma$ is a tree.

Finally, because $F \cap G$ is a connected path, its length
is simply $|F \cap G|$.
\end{proof}

\begin{lemma}\label{leafleaf}
Let $x_{i_1},\ldots,x_{i_t}$ be a path of length $t$ in $\Gamma$
such that $x_{i_t}$ is a leaf of $\Gamma$.
Then the facet $\{x_{i_1},\ldots,x_{i_t}\}$  is a leaf in the simplicial
complex $\Delta_t(\Gamma)$.
\end{lemma}

\begin{proof}
Let $F$ denote both the path and the facet in $\Delta_t(\Gamma)$.
Our proof depends upon whether or not $x_{i_1}$ is the root of
the tree $\Gamma$.  Note that there is nothing to prove
if $F$ is the only path of length $t$ of $\Gamma$.

Suppose that $x_{i_1}$ is not the root.
There then exists another vertex $x_{i_0}$ such that $(x_{i_0},x_{i_1})$
is a directed edge in $\Gamma$.  But this means that
$x_{i_0},\ldots,x_{i_{t-1}}$ is also a path of length $t$ in $\Gamma$.
In other words, $G = \{x_{i_0},\ldots,x_{i_{t-1}}\} \in \Delta_t(\Gamma)$.
But then for any $F \neq F' \in \Delta_t(G)$,
$F' \cap F \subseteq F \cap G = \{x_{i_1},\ldots,x_{i_{t-1}}\}$ because
there is only one path of length $t$ in $\Gamma$ that ends
at $x_{i_t}$.  Thus, $F$ is a leaf of $\Delta_t(\Gamma)$.

Suppose now that $x_{i_1}$ is the root of the tree.  Suppose that
\[\bigcup_{
\begin{array}{c}
G \in \Delta_t(\Gamma) \\
\mbox{$G$ is a facet, $G \neq F$}
\end{array}} (G \cap F)= \{x_{i_1},\ldots,x_{i_m} \}. \]
Note that because there are at least two paths of length $t$ in
$\Gamma$, there are at least two paths that start at the
root $x_{i_1}$.  Thus, $x_{i_1}$ is indeed an element of
the above set, i.e., the union is non-empty.
Let $x_{i_m}$ be the vertex with largest index
in the above union.  So, there exists a facet $H \in \Delta_t(\Gamma)$
that contains $x_{i_m}$.  Furthermore, $x_{i_m} \neq x_{i_t}$
because otherwise there are two different paths that start at the
root and end at
the leaf $x_{i_t}$, contradicting the fact that $\Gamma$
is a tree.

The path $H$ must have the form $H =
x_{i_j},\ldots, x_{i_m},x_{r_{m+1}},\ldots,x_{r_k}$ with  $1\leq j\leq m <
t$ and $k = t+j-1$ since $H$ is a path of length $t$.
Note that $x_{i_j},\ldots,x_{i_m}$ must also be in $F$.  If not,
we could find two distinct paths from the root to $x_{i_m}$.
But this means that $x_{i_1},\ldots,x_{i_j},\ldots,x_{r_k}$
is a path of length $t + j-1$.  So there
is a path of length $t$ from $x_{i_1}$ to $x_{r_t}$, whence
$G = \{x_{i_1},\ldots,x_{i_m},x_{r_{m+1}},\ldots,x_{r_t}\}
\in \Delta_t(\Gamma)$.   It is now straightforward
to verify that $F \cap F' \subseteq F \cap G = \{x_{i_1},\ldots,x_{i_m}\}$
for all
$F' \neq F$ in $\Delta_t(\Gamma)$, and so $F$ is a leaf.
\end{proof}

We have now come to the main theorem of this section.

\begin{theorem}\label{simplicialtree}
Let $\Gamma$ be a tree and $t \geq 2$.
Then $\Delta_t(\Gamma)$ is a simplicial tree.
\end{theorem}
\begin{proof}
We proceed by induction on $\height(\Gamma)$.  Note that
if $\height(\Gamma) < t-1$, then
there are no paths of length $t$ in $\Gamma$, and
so $\Delta_t(\Gamma) = \{\}$.  By convention, this is a tree.
When  $\height(\Gamma) = t-1$, then all the paths of length $t$
must contain both the root and a leaf of $\Gamma$.
Then by Lemma \ref{leafleaf}, all the paths correspond to
leaves of the simplicial complex $\Delta_t(\Gamma)$.  The result
now follows.

We now suppose that $\height(\Gamma) = r > t-1$ and that
the result holds for all trees of height $r'$ with $0 < r' < r$.
Let $w_1,\ldots,w_m$ denote all the leaves of $\Gamma$
whose level equals $r$.  Let $\Gamma' = \Gamma \setminus \{w_1,\ldots,w_m\}$
denote the tree obtained by removing the leaves $w_1,\ldots,w_m$
and their corresponding incident edges.  Note that
$\height(\Gamma') = r-1$.

Let $\langle G_1,\ldots,G_k \rangle$ be any subcomplex of $\Delta_t(\Gamma)$.
If for each $i =1,\ldots,k$, we have $G_i \cap \{w_1,\ldots,w_m\} = \emptyset$,
then each $G_i$ is also a facet of $\Delta_t(\Gamma')$, i.e.,
$\langle G_1,\ldots,G_k\rangle$ is also a subcomplex of $\Delta_t(\Gamma')$.
By induction, this subcomplex contains a leaf.

So, suppose that after relabeling $G_1 \cap \{w_1,\ldots,w_m\} \neq
\emptyset$.  We continue to relabel so that
\[|G_1 \cap G_2 | \leq |G_1 \cap G_3| \leq \cdots \leq |G_1 \cap G_k|. \]
The conclusion will follow once we show that $G_1 \cap G_i \subseteq
G_1 \cap G_k$ for $i=2,\ldots,k$, i.e., $G_1$ is a leaf.  Note that if
$G_1 \cap G_i = \emptyset$, then the conclusion
$G_1 \cap G_i \subseteq G_1 \cap G_k$ is immediate.  So,
suppose that $G_1 \cap G_i \neq \emptyset$.  To simplify
our notation, write $G_1 = \{y_1,\ldots,y_t\}$ where
$y_1,\ldots,y_t$ is a path of length of $t$ that ends at $y_t$
and $y_t \in \{w_1,\ldots,w_m\}$.

By Lemma \ref{intersectionoftwopaths}, $G_1 \cap G_i$ is a connected
path of length $|G_1 \cap G_i|$ that starts at the starting
vertex of $G_1$ or $G_i$.  In particular,
it starts at the vertex with the largest
level.  This vertex must
by $y_1$.  Indeed, if the path $G_i$ starts at a vertex of larger level,
say $G_i = \{v_1,\ldots,v_t\}$ with $\level(v_1) > \level(y_1)$,
then
\begin{eqnarray*}
r = \height(\Gamma) \geq \level(v_t) &= &\level(v_1) + (t-1) \\
& > & \level(y_1) + (t-1) = \level(y_t) = r,
\end{eqnarray*}
thus giving the desired contradiction.

Since $G_1 \cap G_i$ is a connected path of length $e = |G_1 \cap G_i|$,
we must have that $G_1 \cap G_i  = \{y_1,y_2,\ldots,y_e\}$.  This
argument holds for any $G_i$ with $G_1 \cap G_i \neq \emptyset$.
In particular,
\[G_1 \cap G_i = \{y_1,y_2,\ldots,y_{|G_1 \cap G_i|}\} \subseteq
\{y_1,y_2,\ldots,y_{|G_1 \cap G_k|}\} = G_1 \cap G_k,\]
as desired.
\end{proof}

\begin{remark}
Theorem \ref{simplicialtree} provides a means to construct
simplicial trees starting from any tree $\Gamma$.
We point out that
it is not true that every simplicial tree $\Delta$ equals
$\Delta_t(\Gamma)$ for some tree $\Gamma$.  A simplicial
tree need not be {\bf pure}, that is, all the facets
have the same cardinality.  On the other hand, the simplicial
trees of the form $\Delta_t(\Gamma)$ are all pure simplicial
complexes.
\end{remark}

Because $I_t(\Gamma)$ is the facet ideal of $\Delta_t(\Gamma)$,
the following holds:

\begin{corollary}\label{facetidealtree}
Let $\Gamma$ be a tree and $t \geq 2$.  Then $I_t(\Gamma)$ is the
facet ideal of a simplicial tree.
\end{corollary}

There are a number of results in the literature about simplicial trees that
can now be applied to $I_t(\Gamma)$.   We end this section
with one such result about sequentially Cohen-Macaulay modules.

\begin{definition}
If $R = k[x_1,\ldots,x_n]$, then
a graded $R$-module $M$ is called {\bf sequentially Cohen-Macaulay} (over $k$)
if there exists a finite filtration of graded $R$-modules
\[ 0 = M_0 \subset M_1 \subset \cdots \subset M_r = M \]
such that each $M_i/M_{i-1}$ for $i = 1,\ldots,r$
is Cohen-Macaulay, and the Krull dimensions of the
quotients are increasing:
$\dim (M_1/M_0) < \dim (M_2/M_1) < \cdots < \dim (M_r/M_{r-1}).$
\end{definition}

Faridi \cite{F2} showed that simplicial trees give
rise to sequentially Cohen-Macaulay modules:

\begin{theorem}[{\cite[Corollary 5.6]{F1}}]\label{simtree SCM}
If $I = I(\Delta)$ is
the facet ideal of a simplicial tree $\Delta$, then
$R/I$ is
sequentially Cohen-Macaulay.
\end {theorem}

Combining Corollary \ref{facetidealtree} and Theorem \ref{simtree SCM},
we then get:

\begin{corollary}\label{scm}
Let $\Gamma$ be a tree and $t \geq 2$.  Then
$R/I_t(\Gamma)$ is sequentially Cohen-Macaulay.
\end{corollary}


\section{Some homological invariants of path ideals of trees}

In this section we investigate some of the invariants of
the ideal $I_t(\Gamma)$ that are encoded into the minimal free
graded resolution of $I_t(\Gamma)$.  A highlight of
this section is a recursive formula for computing the projective
dimension of $I_t(\Gamma)$ if the simplicial complex $\Delta_t(\Gamma)$
is properly-connected.

Let $M$ be any graded $R$-module where $R = k[x_1,\ldots,x_n]$.  Associated
to $M$ is a {\bf minimal free graded resolution} of the form
\[
0 \rightarrow \bigoplus_j R(-j)^{\beta_{p,j}(M)}
\rightarrow \bigoplus_j R(-j)^{\beta_{p-1,j}(M)}
\rightarrow \cdots
\rightarrow  \bigoplus_j R(-j)^{\beta_{0,j}(M)}
\rightarrow M \rightarrow 0
\]
where the maps are exact, $p \leq n$, and $R(-j)$ is the
$R$-module obtained by shifting the degrees of $R$ by $j$.
The number $\beta_{i,j}(M)$, the $(i,j)$-th
{\bf graded Betti number} of $M$, is an invariant of $M$ that
equals the number of minimal generators of degree $j$
in the $i$th syzygy module.  The {\bf projective dimension}
of $M$, denoted $\pd(M)$, is equal to $p$, the minimal
length of all free resolutions of $M$.

For a general monomial ideal $I$, the numbers $\beta_{i,j}(I)$
may depend upon the characteristic of the field $k$.  However,
this is not the case for the ideals of the form $I_t(\Gamma).$

\begin{theorem}\label{characteristic}
Let $\Gamma$ be a tree and $t \geq 2$.  Then the graded Betti
numbers of $I_t(\Gamma)$ are independent of
the characteristic of $k$.
\end{theorem}

\begin{proof}
\cite[Theorem 5.8]{HVT2} proved the existence of a recursive formula
(for brevity, we do not reproduce this formula)
to compute the graded Betti numbers of the facet ideal of
simplicial tree $\Delta$
in terms of the graded Betti numbers of subcomplexes
of $\Delta$ that are also trees.  This formula
is independent of the characteristic; indeed,
the recursive nature of the formula allows
us to reduce to the case that the facet ideal has single
generator, and the resolution of this ideal does not depend
upon the characteristic.  The conclusion
follows from Corollary \ref{facetidealtree}.
\end{proof}

We now derive a recursive formula for $\pd(I_t(\Gamma))$
in a special case
using some tools developed in \cite{FHVT1}.  We let $G(I)$
denote the unique set of minimal generators of a monomial ideal $I$.

\begin{definition}
Let $I$ be a monomial ideal, and suppose that there exists
monomial ideals $J$ and $K$ such that $G(I)$ is
the disjoint union of $G(J)$ and $G(K)$.  Then
$I = J+K$ is a {\bf Betti splitting} if
\[\beta_{i,j}(I) = \beta_{i,j}(J) + \beta_{i,j}(K) + \beta_{i-1,j}(J\cap K)
~~\mbox{for all $i,j \geq 0$}.\]
\end{definition}

In the paper \cite{FHVT1}, the authors describe a number of
sufficient conditions for an ideal $I$ to have a Betti splitting.  We
shall require the following such condition:

\begin{theorem}[{\cite[Corollary 2.7]{FHVT1}}]\label{bettisplit}
Suppose that $I = J+K$ where $G(J)$ contains all the generators
of $I$ divisible by the variable $x_i$ and $G(K)$ is a nonempty set
containing the remaining
generators of $I$.  If $J$ has a linear resolution, then $I = J+K$
is a Betti splitting.
\end{theorem}

\begin{corollary}\label{bettisplittree}
Let $\Gamma$ be a tree with $t \geq 2$.  Suppose that
$w = x_{i_1}\cdots x_{i_t}$ is a generator of $I_t(\Gamma)$,
and furthermore, suppose that $x_{i_t}$ is a leaf of $\Gamma$.
Then
\[I_t(\Gamma) = (w) + I_t(\Gamma')  ~~\mbox{with $\Gamma' = \Gamma \setminus
\{x_{i_t}\}$}\]
is a Betti splitting of
$I_t(\Gamma)$.  Consequently,
\[\pd(I_t(\Gamma)) = \max\{\pd(I_t(\Gamma')),\pd((w) \cap I_t(\Gamma'))+1\} .\]
\end{corollary}

\begin{proof}
The ideal $J = (w)$ is generated by all the generators of $I_t(\Gamma)$
divisible by $x_{i_t}$.  There is only one since $x_{i_t}$ is a leaf.
The remaining generators of $I_t(\Gamma)$ generate $I_t(\Gamma')$.
Because $J$ has a linear resolution,  Theorem \ref{bettisplit} applies.
The statement about the projective dimension follows from
Corollary 2.2 of \cite{FHVT1}, and the fact that $\pd((w)) =0$.
\end{proof}

In general, the projective dimension of the ideal $(w) \cap I_t(\Gamma')$
may be hard to compute.  However, if we assume that $\Delta_t(\Gamma)$
has an extra property, we can describe a recursive formula.  We introduce
the notion of a properly-connected simplicial complex, which was first
defined for hypergraphs by H\`a and the second author \cite{HVT1}.

\begin{definition} Let $\Delta$ be a pure simplicial complex
where every facet has dimension $d$. A {\bf chain of length $n$} in
$\Delta$ is a sequence of facets $(F_0,\ldots,F_n)$ such that
\begin{enumerate}
\item[$(1)$]  $F_0,\ldots,F_n$ are all distinct facets of $\Delta$, and
\item[$(2)$]  $F_i \cap F_{i+1} \neq \emptyset$ for
$i=0,\ldots,n-1$.
\end{enumerate}
Two facets are $F$ and $F'$ are {\bf connected} if there
exists a chain $(F_0,\ldots,F_n)$ where $F = F_0$ and $F' = F_n$.
The chain connecting $F$ to $F'$ is a {\bf proper chain} if
$|F_i \cap F_{i+1}| = |F_{i+1}|-1$ for all $i = 0,\ldots,n-1$.
The (proper) chain is an {\bf (proper) irredundant chain} of length $n$
if no proper subsequence is a (proper) chain from $F$ to $F'$.
The {\bf distance} between two facets $F$ and $F'$ in $\Delta$
is then given by
\[\dist_{\Delta}(F,F') = \min\{\ell ~|~ (F= F_0,\ldots,F_{\ell}=F') ~\mbox{is a
proper irredundant chain}\};\] if no such chain exists, then set
$\dist_{\Delta}(F,F') = \infty$. We say that $\Delta$ is {\bf
properly-connected} if 
\[\dist_{\Delta}(F,F') = (d+1) - |F\cap F'|\]
for any two facets $F,F' \in \Delta$ with the
property that $F \cap F' \neq \emptyset$.
Otherwise, we say $\Delta$
is {\bf not properly-connected}.
\end{definition}

We need one more preparatory lemma and some notation before coming to
the main result of this section.
Let $\Gamma$ be a tree and $t \geq 2$.
Suppose that $F = \{x_{i_1},\ldots,x_{i_t}\} \in \Delta_t(\Gamma)$
with $x_{i_t}$ a leaf of $\Gamma$.  Let $Y$ denote the set
\[Y := \bigcup_{\begin{array}{c}
G \in \Delta_t(\Gamma) \\
|G \cap F| = t-1
\end{array}} (G \setminus F) = \{y_1,\ldots,y_a\}
\]
and let $Z := Y \cup F = \{y_1,\ldots,y_a,x_{i_1},\ldots,x_{i_t}\}$.

\begin{lemma}  \label{intersectionlemma}
Let $\Gamma$ be a tree with $t \geq 2$, and suppose that
$\Delta_t(\Gamma)$ is properly connected.
Suppose that $w = x_{i_1}\cdots x_{i_t} \in I_t(\Gamma)$ is
a monomial generator of $I_t(\Gamma)$ with $x_{i_t}$ a leaf.
Set $\Gamma' = \Gamma \setminus \{x_{i_t}\}$ and
$\Gamma'' = \Gamma \setminus Z$, with $Z$ as above.  Then
\[I_t(\Gamma') \cap (w) = w(y_1,\ldots,y_a) + wI_t(\Gamma'')\]
with $Y = \{y_1,\ldots,y_a\}$ as above.
\end{lemma}

\begin{proof}
This lemma is now a special case of \cite[Corollary 4.14]{HVT1} once one
translates the statements from hypergraphs to the simplicial complex
case.
\end{proof}

We can now prove the following formula for the projective
dimension.  Recall that $\pd(R/I) = \pd(I)+1$, and we
adopt the convention that when $I = (0)$, then $\pd(I) = -1$.

\begin{theorem}\label{generalprojdim}
Let $\Gamma$ be a tree with $t \geq 2$, and furthermore,
suppose that $\Delta_t(\Gamma)$ is properly-connected.  Let
$w = x_{i_1}\cdots x_{i_t}$ be a generator of $I_t(\Gamma)$
with $x_{i_t}$ a leaf.    Then
\[\pd(I_t(\Gamma)) = \max\{\pd(I_t(\Gamma')),\pd(I_t(\Gamma''))+a+1\}\]
where $\Gamma' = \Gamma \setminus \{x_{i_t}\}$, $\Gamma'' = \Gamma \setminus
Z$, and $a = |Y|$, where $Z$ and $Y$ are defined as above.
\end{theorem}

\begin{proof}  It suffices by Corollary \ref{bettisplittree}
to show that  $\pd((I_t(\Gamma')\cap (w))
= \pd(I_t(\Gamma''))+a$.
By Lemma \ref{intersectionlemma} we know that
\[I_t(\Gamma') \cap (w) =
w\left[(y_1,\ldots,y_a)+ I_t(\Gamma'')\right].\]
Because none of the variables that divide $w$
divide any  generator of $(y_1,\ldots,y_a)+ I_t(\Gamma'')$,
\[\pd(I_t(\Gamma') \cap (w)) =
\pd((y_1,\ldots,y_a)+ I_t(\Gamma'')).\]
Since no generator of $I_t(\Gamma'')$ is divisible by
any element of $\{y_1,\ldots,y_a\}$, we have the isomorphism
\[k[y_1,\ldots,y_a,w_1,\ldots,w_b]/((y_1,\ldots,y_a)+ I_t(\Gamma''))
\cong k[y_1,\ldots,y_a]/(y_1,\ldots,y_a) \otimes k[w_1,\ldots,w_b]/I_t(\Gamma'')\]
where $\{w_1,\ldots,w_b\}  = V \setminus Y$.  Because $y_1,\ldots,y_a$
form a regular sequence on $k[y_1,\ldots,y_a]$, we have
$\pd(k[y_1,\ldots,y_a]/(y_1,\ldots,y_a)) = a$.

Putting together these pieces, we now have
\begin{eqnarray*}
\pd(R/((y_1,\ldots,y_a)+I_t(\Gamma''))) & = & \pd(k[y_1,\ldots,y_a]/(y_1,\ldots,y_a) \otimes k[w_1,\ldots,w_b]/I_t(\Gamma'')) \\
& = & \pd(k[y_1,\ldots,y_a]/(y_1,\ldots,y_a)) + \pd(k[w_1,\ldots,w_b]/I_t(\Gamma'')) \\
& = & a + \pd(I_t(\Gamma'')) +1.
\end{eqnarray*}
It thus follows that
$$\pd((y_1,\ldots,y_a)+I_t(\Gamma'')) =
\pd(R/((y_1,\ldots,y_a)+I_t(\Gamma''))) - 1 = a+ \pd(I_t(\Gamma'')),$$
thus finishing the proof.
\end{proof}

\begin{remark}
When $t = 2$, then $I_t(\Gamma)$ is the edge ideal of $\Gamma$,
and the facets of $\Delta_t(\Gamma)$ are simply the edges of $\Gamma$.
As remarked in \cite[Example 4.4]{HVT1}, every simple graph is
properly-connected, so the formula in Theorem \ref{generalprojdim}
holds for all $\Gamma$, and one can show that our result
is the same as the formula for the projective dimension of trees
as first found in \cite{JK1}.
\end{remark}

\section{Case Study:  the path ideal of the line graph}

In this section we study $I_t(\Gamma)$ in the
special case that $\Gamma$ is the line graph $L_n$.
The {\bf line graph} $L_n$ is a tree with vertex set $V = \{x_1,\ldots,x_n\}$
and directed edges $e_j = (x_j,x_{j+1})$ for $j=1,\ldots,n-1$.
Thus, the graph $L_n$ has the form
\[
\begin{picture}(200,20)(0,0)
\put(0,10){\circle*{5}}
\put(0,15){$x_1$}
\put(0,10){\line(1,0){100}}
\put(0,10){\vector(1,0){20}}
\put(0,10){\vector(1,0){60}}
\put(0,10){\vector(1,0){100}}
\put(105,10){\line(1,0){10}}
\put(120,10){\line(1,0){10}}
\put(135,10){\line(1,0){10}}
\put(150,10){\line(1,0){50}}
\put(160,10){\vector(1,0){20}}
\put(40,10){\circle*{5}}
\put(40,15){$x_2$}
\put(80,10){\circle*{5}}
\put(80,15){$x_3$}
\put(160,10){\circle*{5}}
\put(160,15){$x_{n-1}$}
\put(200,10){\circle*{5}}
\put(200,15){$x_n$}
\end{picture}
\]
In this section we give an exact formula for
$\pd(I_t(L_n))$ in terms of $t$ and $n$,
we use this formula to compute the arithmetical rank
of $I_t(L_n)$ for some $t$ and $n$, and we show
that this formula provides a lower bound on $\pd(I_t(\Gamma))$
for any rooted tree $\Gamma$ of height $n-1$.

\subsection{Projective dimension}
It is not hard to verify that for all integers $n \geq t \geq 2$,
the simplicial complex $\Delta_t(L_n)$ is properly-connected.
We can then use Theorem \ref{generalprojdim}
to find an exact formula for the projective dimension.

\begin{theorem}\label{maintheorem1}
Fix integers $n \geq t \geq 2$.  Then
\[\operatorname{pd}(R/I_t(L_n))=
\left\{
\begin{array}{ll}
\displaystyle{\frac{2(n-d)}{t+1}} & \hbox{if $n\equiv d \pmod{(t+1)}$ with $0 \leq d \leq t-1$} \\
\displaystyle{\frac{2n - (t-1)}{t+1}} & \hbox{if $n\equiv t \pmod{(t+1)}$.}
\end{array}
\right.
\]
\end{theorem}

\begin{proof}
To prove that the formula holds for all $(n,t)$ with $n \geq t \geq 2$,
we will use induction on the tuple $(n,t)$.
For any $t$, if $n=t$, then $I_t(L_n) = (x_1\cdots x_n)$ is a principal
ideal, and thus $\pd(R/I_t(L_m)) = 1$.  This agrees with the formula
in statement.

We now prove the statement for all $(n,t)$ with $n>t$, and by induction,
we can assume that the statement holds for all $(k,t)$ with $n > k \geq t$.

The monomial $x_{n-t+1}x_{n-t+2}\cdots x_n$ is a generator
of $I_t(L_n)$, and furthermore, $x_n$ is a leaf of $L_n$.  There
is one path in $L_n$ of length $t$ that intersects the
path $x_{n-t+1},\ldots,x_n$ at $t-1$ vertices, namely
the path $x_{n-t},\ldots,x_{t-1}$.  Thus, in the notation
of Theorem \ref{generalprojdim}, we have $Y = \{x_{n-t}\}$
and $Z = \{x_{n-t},\ldots,x_n\}$.  But then
$L_n \setminus \{x_n\} = L_{n-1}$ and $L_n \setminus Z = L_{n-t-1}$.
Applying Theorem \ref{generalprojdim} gives us the formula
\[\pd(R/I_t(L_n)) = \max\{ \pd(R/I_t(L_{n-1})), \pd(R/I_t(L_{n-t-1}))+2\},\]
which we shall use in the induction step.

We consider two main cases, depending upon the value of $n \pmod{(t+1)}$.

\noindent
{\bf Case 1.}  $n \equiv d \pmod{(t+1)}$ with $0 \leq d \leq t-1$.

We further subdivide this case into two subcases:

\noindent
{\bf Case 1(a).} $d= 0$.

\noindent
When $d = 0$, then $n-1 \equiv -1 \equiv t \pmod{(t+1)}$.  By induction
\[\pd(R/I_t(L_{n-1})) = \frac{2(n-1)-(t-1)}{t+1} = \frac{2n-t-1}{t+1}.\]
Also, since $d=0$, $n-t-1 \equiv d \equiv 0 \pmod{(t+1)}$.  Either
$n-t-1 \geq t$, in which case induction gives
\[\pd(R/I_t(L_{n-t-1})) = \frac{2(n-t-1-0)}{t+1} = \frac{2n}{t+1} -2,\]
or $n-t-1 < t$ which implies $\pd(R/I_t(L_{n-t-1}))=0$.
In the first case
\[\pd(R/I_t(L_n)) = \max\left\{\frac{2n-t-1}{t+1},\frac{2n}{t+1}-2+2\right\}
= \frac{2n}{t+1} =
\frac{2(n-d)}{t+1}\]
where the last equality comes from the fact that $d=0$.
In the second situation, because
$t < n < 2t+1$ and $n \equiv 0 \pmod{(t+1)}$, this actually forces
$n = t+1$.  But then
\[\pd(R/I_t(L_n)) = \max\left\{\frac{t+1}{t+1},2\right\} = 2 =
\frac{2(n-d)}{t+1}\]
where the last equality again comes from the fact that $d=0$
and $n=t+1$.
So the formula agrees in these cases.

\noindent
{\bf Case 1(b).} $d \neq 0$.

\noindent
Because $d \neq 0$, then $n-1 \equiv d-1 \pmod{(t+1)}$
with $0 \leq d-1 \leq t-1$.  By
induction
\[\pd(R/I_t(L_{n-1}))=\frac{2(n-1-(d-1))}{t+1}=\frac{2(n-d)}{t+1}.\]
Also, $n-t-1\equiv d \pmod{(t+1)}$.  Now, either $n-t-1 < t$,
in which case $\pd(R/I_t(L_{n-t-1}))=0$, or $n-t-1 \geq t$,
and thus, by induction
\[\pd(R/I_t(L_{n-t-1}))=\frac{2(n-t-1-d)}{t+1}=\frac{2(n-d)-2(t+1)}{t+1} = \frac{2(n-d)}{t+1}-2.\]
If $\pd(R/I_t(L_{n-t-1}))=0$, then
\[\pd(R/I_t(L_n))=\max\left\{\frac{2(n-d)}{t+1},2\right\}=\frac{2(n-d)}{t+1} = 2\]
where the last equality comes from the fact that when $t < n < 2t+1$
and $n \equiv d \pmod{(t+1)}$, then $n = t+1+d$.
In the second situation,
\[\pd(R/I_t(L_n))=\max\left\{\frac{2(n-d)}{t+1},\frac{2(n-d)}{t+1}\right\}=\frac{2(n-d)}{t+1}.\]
Again, the formula agrees with the statement.

\noindent
{\bf Case 2.} $n \equiv t \pmod{(t+1)}$.

\noindent
In this case $n-1 \equiv t-1 \pmod{(t+1)}$, thus by induction
\[\pd(R/I_t(L_{n-1}))=\frac{2(n-1-(t-1))}{t+1}=\frac{2(n-t)}{t+1}.\]
Note that $n-t-1 \not< t$, because if it was,
this would mean that $t < n < 2t+1$ and $n \equiv t \pmod{(t+1)}$.
This, in turn, would imply that $(n-t) = k(t+1)$ for some $k$,
i.e., $n = (k+1)t+k$.  But there is no integer $k$ that gives
$t < (k+1)t+k < 2t+1$.  So $\pd(R/I_t(L_{n-t-1})) \neq 0$,
and moreover, induction implies,
\[\pd(R/I_t(L_{n-t-1}))=\frac{2(n-t-1)-(t-1))}{t+1}=\frac{2n-(t-1)}{t+1} -2.\]
So
\[\pd(R/I_t(L_n))=\max\left\{\frac{2(n-t)}{t+1},\frac{2n-(t-1)}{t+1} \right\}=
\frac{2n-(t-1)}{t+1}.\]
Again, this agrees with formula in the statement.

These three cases now complete the proof.
\end{proof}

Jacques \cite{J1} gave a formula for the projective
dimension of the edge ideal of the graph $L_n$ for all $n$.
This formula
is now a corollary of
Theorem \ref{maintheorem1} by specializing to the case that
$t=2$.

\begin{corollary}[{\cite[Corollary 7.7.35]{J1}}] Let $n \geq 2$ be an integer, and let $L_n$ be the line graph.
 Then
\[\operatorname{pd}(R/I_2(L_n))=
\left\{
\begin{array}{ll}
\displaystyle{\frac{2(n-d)}{3}} & \hbox{if $n\equiv d \pmod{3}$ with $0 \leq d \leq 1$} \\
\displaystyle{\frac{2n-1}{3}} & \hbox{if $n\equiv 2 \pmod{3}$.}
\end{array}
\right.
\]
\end{corollary}

\subsection{Arithmetical rank}  The {\bf arithmetical rank}
of an ideal $I \subseteq R$, denoted $\ara(I)$, is the least number of
elements of $R$ that generate $I$ up to radical, that is,
\[\ara(I) = \min\left\{s ~\left|~ ~\mbox{there exist $f_1,\ldots,f_s \in R$ such
that $\sqrt{I} = \sqrt{(f_1,\ldots,f_s)}$}\right\}\right..\]
Determining $\ara(I)$ for an arbitrary ideal (even monomial
ideals!) of $R$ remains an open
and difficult problem.  We will compute $\ara(I_t(L_n))$ for some
values of $t$ and $n$.

When $I$ is a square-free monomial ideal, then
Lyubeznik \cite{L1} proved that
\[\pd(R/I) \leq \ara(I).\]
Theorem \ref{maintheorem1} thus gives us a lower bound
on $\ara(I_t(L_n))$ for all $t$ and $n$.
The following result of Schmitt and Vogel is one
of the few known techniques available for computing
an upper bound on the arithmetical rank of an ideal.

\begin{theorem}[\cite{SV1}]\label{SchmittVogel}
Let $P$ be a finite subset of elements of $R$. Let $P_0,\ldots,P_r$
be subsets of $P$ such that
\begin{enumerate}
\item $\bigcup_{i=0}^{r}P_i=P$;
\item $P_0$ has exactly one element;
\item if $p$ and $p^{\prime}$ are
  different elements of $P_i$ ($0<i\leq
  r$), then there is an integer $i^{\prime}$
  with $0\leq i^{\prime}<i$ and an
  element in $P_{i^{\prime}}$ which divides
  $pp^{\prime}$.
\end{enumerate}
Set $q_i=\sum_{p\in P_i}p^{e(p)}$, where $e(p)\geq 1$ are
arbitrary integers, and write $(P)$ for the ideal of $R$
generated by the elements of $P$. Then
\[\sqrt{(P)}=\sqrt{(q_0,\ldots,q_r)}.\]
In particular, $\ara((P)) \leq r+1$.
\end{theorem}

One strategy that has been used
to compute $\ara(I)$ for a monomial
ideal is to use Theorem \ref{SchmittVogel} to partition
the generators of $I$  into $\pd(R/I)$ subsets that satisfy
the conditions $(1)-(3)$ of the theorem.  Theorem \ref{SchmittVogel}
and Lyubeznik's lower bound then shows that $\pd(R/I) = \ara(I)$.
For example, this strategy is used by Barile \cite{B1}
to compute $\ara(I_2(L_n))$ for all $n$:

\begin{theorem}[\cite{B1}]  Let $L_n$ be the line graph with $n \geq 2$.  Then
\[\ara(I_2(L_n)) = \pd(R/I_2(L_n)).\]
\end{theorem}

We also employ this strategy to prove the following theorem.

\begin{theorem} \label{t3case}
Let $L_n$ be the line graph with $n \geq 3$.
Let $t = 3$ and suppose that $n \not\equiv 2 \pmod{4}$.  Then
\[\ara(I_3(L_n))=\pd(R/I_3(L_n)).\]
\end{theorem}

\begin{proof}
We shall break this proof into three cases.

\noindent
{\bf Case 1.} $n \equiv 0 \pmod{4}$.

\noindent
In this case $n = 4k$ for some integer $k$.
By Theorem \ref{maintheorem1} we have
$\pd(R/I_3(L_n))=\frac{2(4k)}{4}=2k$.  On the other hand,
the number of generators of $I_3(L_n)$ is $n-3+1 = 4k-3+1=4k-2$.
We partition the $4k-2$ generators into the following $2k$ sets:
\begin{eqnarray*}
P_0 &=& \{x_2x_3x_4\} \\
P_i&=&\{x_{2i-1}x_{2i}x_{2i+1},x_{2i+2}x_{2i+3}x_{2i+4}\}
~~~\mbox{for  $i=1,\ldots,2k-2$} \\
P_{2k-1}&=&\{x_{4k-3}x_{4k-2}x_{4k-1}\}.
\end{eqnarray*}
Because this partition of the generators of $I_3(L_n)$
satisfies the conditions of
Theorem \ref{SchmittVogel}, we have $2k = \pd(R/I_3(L_n))
\leq \ara(I_3(L_n)) \leq 2k$.

\noindent
{\bf Case 2.} $n \equiv 1 \pmod{4}$.

\noindent
In this situation, $n = 4k+1$.  By Theorem \ref{maintheorem1},
we have $\pd(R/I_3(L_n))=\frac{2(4k+1-1)}{4}=2k$.
The number of generators of $I_3(L_n)$ is $n-3+1 = 4k+1-3+1=4k-1$.
We partition the $4k-1$ generators into the following $2k$ sets:
\begin{eqnarray*}
P_0 &=& \{x_2x_3x_4\}\\
P_i &=& \{x_{2i-1}x_{2i}x_{2i+1},x_{2i+2}x_{2i+3}x_{2i+4}\},
~~\mbox{for $i=1,\ldots,2k-2$} \\
P_{2k-1}&=& \{x_{4k-3}x_{4k-2}x_{4k-1},x_{4k-1}x_{4k}x_{4k+1}\}.
\end{eqnarray*}
This partition satisfies the hypotheses of Theorem \ref{SchmittVogel}.
The conclusion now follows.

\noindent
{\bf Case 3.}  $n \equiv 3 \pmod{4}$.

\noindent
We can write $n$ as $n = 4k+3$.  The projective
dimension of $R/I_3(L_n)$ is then
$\pd(R/I_3(L_n))=\frac{2(4k+3) - 2}{4}=\frac{8k+4}{4}=2k+1$
by Theorem \ref{maintheorem1}.  The ideal $I_3(L_n)$
has $n-3+1 = 4k+3-3+1 = 4k+1$ minimal generators.
We partition these generators into the following $2k+1$  sets:
\begin{eqnarray*}
P_0 &=&\{x_2x_3x_4\} \\
P_i &=& \{x_{2i-1}x_{2i}x_{2i+1},x_{2i+2}x_{2i+3}x_{2i+4}\}
~~\mbox{for $i=1,\ldots,2k-1$}\\
P_{2k} &=&\{x_{4k-1}x_{4k}x_{4k+1},x_{4k+1}x_{4k+2}x_{4k+3}\}.
\end{eqnarray*}
If we now
apply the theorem of Schmitt and Vogel, the desired
conclusion follows.
\end{proof}

When $n \equiv 2 \pmod{4}$, we were not able to compute
the arithmetical rank of the ideal $I_3(L_n)$.  As we will
show below, there is no way to partition the {\it generators}
of $I_3(L_n)$ into the $p = \pd(R/I_3(L_n))$ sets that satisfy
the hypotheses of Theorem \ref{SchmittVogel}.  However,
there may still exist $p$ elements $f_1,\ldots,f_p$ of $R$
such that $I_3(L_n) = \sqrt{(f_1,\ldots,f_p)}$.

\begin{definition}
Let $I$ be a monomial ideal of $R$ and let $G(I)$
denote its minimal generators.
A partition $P_0,\ldots,P_r$ of $G(I)$, that
is, $G(I)$ is the disjoint union of $P_0,\ldots,P_r$,
is a {\bf good partition of $I$}
if $P_0,\ldots,P_r$ also satisfy $(1)-(3)$ of Theorem \ref{SchmittVogel}
and $r+1 = \pd(R/I)$.
\end{definition}

When $I = I_t(L_n)$ we can make some general statements
about a good partition of $I$.  In the following,
let $m_i = x_{i}\cdots x_{i+t-1}$, and thus
 $I_t(L_n)  = (m_1,m_2,\ldots,m_{n-t+1})$.

\begin{lemma}
Let $P_0,\ldots,P_r$ be a good partition of $I_t(L_n)$.
\begin{enumerate}
\item If $m_i, m_j \in P_{\ell}$ with $i < j$, then
$i+1 < j \leq i+t$.
\item $|P_{\ell}| \leq \left\lfloor \frac{t}{2} \right\rfloor +1$
for all $\ell = 1,\ldots,r$.
\end{enumerate}
\end{lemma}
\begin{proof}
(1) If $m_i$ and $m_{i+1}$ are in $P_{\ell}$, then
 the only generators of  $I_t(L_n)$ that divide $m_im_{i+1}
= x_i(x_{i+1}\cdots x_{i+t-1})^2x_{i+t}$ are $m_i$ and $m_{i+1}$.  Thus,
since $P_0,\ldots,P_r$ is a good partition, and hence a partition
of $G(I)$,  condition (3)
of Theorem \ref{SchmittVogel} is violated.  Similarly,
if $m_i$ and $m_j \in P_{\ell}$ with $j > i+t$,
the only generators of $I_t(L_n)$ that divide $m_im_j$ are
$m_i$ and $m_j$, and since the $P_i$'s form a partition,
$m_i$ and $m_j$ do not belong to a $P_t$ with $0 \leq t < \ell$.

(2)
Suppose that $m_i \in P_{\ell}$ and $m_i$ has the smallest index.
By (1), $P_{\ell}$ must be a subset of $\{m_i,m_{i+1},\ldots,m_{i+t}\}$.
Again by $(1)$, $P_{\ell}$ cannot contain ``adjacent'' generators,
i.e., $m_j$ and $m_{j+1}$.  So $P_{\ell}$ can contain at most
half of the elements $\{m_i,\ldots,m_{i+t}\}$.  Thus, it can
contain at most $\left\lfloor \frac{t}{2} \right\rfloor +1$ elements.
\end{proof}

\begin{corollary} \label{nogoodpartition}
 Suppose $n \geq t \geq 2$ are integers such that
\[\left(\pd(R/I_t(L_n))-1\right)\left(\left\lfloor \frac{t}{2} \right\rfloor +1 \right)
< n -t.\]
Then there does not exist a good partition of the generators
of $I_t(L_n)$.
\end{corollary}

\begin{proof}
Suppose that $P_0,\ldots,P_r$ was a good partition of the generators
with $r +1 = \pd(R/I_t(L_n))$.  Because $|G(I_t(L_n))| = n-t+1$
\begin{eqnarray*}
1+ n-t &=& |P_0 \cup P_1 \cup \cdots \cup P_r| = |P_0|+\cdots +|P_r| \\
 &\leq  & |P_0| + \left(\left\lfloor\frac{t}{2}\right\rfloor +1\right)
+ \cdots +
\left(\left\lfloor \frac{t}{2} \right\rfloor+1\right) ~~\mbox{by
Lemma \ref{nogoodpartition}}\\
&=& 1 + r\left(\left\lfloor \frac{t}{2} \right\rfloor+1\right) = 1 +
 \left(\pd(I_t(L_n))-1\right)
\left(\left\lfloor \frac{t}{2} \right\rfloor + 1 \right) < 1+n-t.
\end{eqnarray*}
The contradiction is now clear.
\end{proof}

We can now show that the procedure used in Theorem \ref{t3case}
will not work for $n \equiv 2 \pmod{4}$.

\begin{corollary}
Suppose that $n \equiv 2 \pmod{4}$.  Then there
does not exist a good partition of $I_3(L_n)$.
\end{corollary}

\begin{proof}
The hypotheses imply that $n = 4k+2$, and thus by Theorem
\ref{maintheorem1}, we have $\pd(R/I_3(L_n)) = \frac{2(n-2)}{4}
= \frac{2(4k)}{4} = 2k$.  Now use Corollary \ref{nogoodpartition}
and the inequality
\[ (\pd(R/I_3(L_n))-1)\left(\left\lfloor \frac{3}{2} \right\rfloor +1\right)
=(2k-1)(2) = 4k-2 <  4k-1 = n-3.\]
\end{proof}


\subsection{Application}  The path ideal of a line
graph $L_n$ is, admittedly, a very special family of path ideals.
However, we can use the results in
the previous two subsections to provide a lower bound on
the projective dimension and arithmetical rank for the path
ideal of any tree $\Gamma$.

\begin{theorem}\label{lowerbound}
Let $\Gamma$ be a rooted tree with
$\height(\Gamma) = n-1$.  Then for any $n \geq t \geq 2$,
\[\pd(R/I_t(L_n)) \leq \pd(R/I_t(\Gamma)) \leq \ara(I_t(\Gamma)).\]
\end{theorem}

\begin{proof}
As observed in the previous subsection, the second inequality
holds for any square-free monomial ideal.  It therefore
suffices to proof the first inequality.

Since $\Gamma$ is a tree of height $n-1$, there exists
a path of length $n$ in $\Gamma$ from the root to
one of the leaves of $\Gamma$.  Let $W = \{w_1,\ldots,w_n\}$
denote the vertices of this path.  Note that
the induced graph of $\Gamma$ on $W$ is precisely the line graph
$L_n$.  We abuse notation and let $L_n$ denote this induced graph.

Because $I_t(L_n)$, respectively $I_t(\Gamma)$, is a square-free monomial
ideal, by the Stanley-Reisner correspondence, it can
be associated to a simplicial complex $\Lambda_t(L_n)$, respectively
$\Lambda_t(\Gamma)$.
If $\Lambda$ is any simplicial complex on $V$, and if
$Y \subseteq V$, we let $\Lambda_Y$ denote the subsimplicial
complex consisting of all the faces of $\Lambda$ contained
in $Y$.   It can then be verified that for any $Y \subseteq W$,
$\Lambda_t(L_n)_Y = \Lambda_t(\Gamma)_Y$.

Let $\tilde{H}_i(\Lambda,k)$ denote the $i$-th reduced simplicial
homology of a simplicial complex $\Lambda$.  Then,
by Hochster's formula
(see, for example \cite[Theorem 5.5.1]{BH}), for any $i,j \geq 0$, we have
\begin{eqnarray*}
\beta_{i,j}(I_t(\Gamma)) &=& \sum_{Y \subseteq V,~~|Y|=j} \dim_k \tilde{H}_{j-i-2}(\Lambda_t(\Gamma)_Y,k) \\
&\geq&   \sum_{Y \subseteq W \subseteq V,~~|Y|=j} \dim_k \tilde{H}_{j-i-2}(\Lambda_t(\Gamma)_Y,k)
\\
&=& \sum_{Y \subseteq W,~~|Y|=j} \dim_k \tilde{H}_{j-i-2}(\Lambda_t(L_n)_Y,k)
= \beta_{i,j}(I_t(L_n)).
\end{eqnarray*}
Because $\pd(I) = \max\{i ~|~ \beta_{i,j}(I) \neq 0\}$
for any homogeneous ideal $I$, and
since $\beta_{i,j}(I_t(\Gamma)) \geq \beta_{i,j}(I_t(L_n))$ for
all $i,j \geq 0$, the conclusion $\pd(I_t(\Gamma)) \geq \pd(I_t(L_n))$,
or equivalently $\pd(R/I_t(\Gamma)) \geq \pd(R/I_t(L_n))$, now
follows.
\end{proof}

\noindent
{\bf Acknowledgments}  The second author acknowledges the support
of NSERC.   This paper is partially based upon the M.Sc. project
of the first author.  We also thank T.H. H\`a and the referee for
their comments.  We would also like to thank Sara Saeedi for
pointing out that we needed an extra hypothesis in Lemma 3.6.


\end{document}